\documentclass[12pt]{article}
\usepackage{amssymb,amsmath}

\newcommand{\Z}{{\mathbb Z}}
\newcommand{\F}{{\mathbb F}}

\renewcommand{\O}{\mathcal{O}}
\newcommand{\s}{\mathcal{S}}
\newcommand{\f}{\mathcal{F}}

\newcommand{\E}{\mathcal{E}}
\newcommand{\B}{\mathcal{B}}
\newcommand{\A}{\mathcal{A}}
\newcommand{\M}{\mathcal{M}}

\newcommand{\T}{\mathcal{T}}

\newcommand{\n}{\mbox{N}}

\newcommand{\Gammabar}{\overline{\Gamma}}
\newcommand{\gammabar}{\overline{\gamma}}
\renewcommand{\char}{\mbox{char}}
\newcommand{\Gal}{\mathrm{Gal}}
\newcommand{\Aut}{\mathrm{Aut}}
\newcommand{\Hom}{\mathrm{Hom}}

\newcommand{\proof}{\noindent{\em Proof: }}
\newcommand{\qed}{\hspace{\fill}$\square$}
\newcommand{\ra}{\rightarrow}
\newcommand{\lra}{\longrightarrow}
\newcommand{\dst}{\displaystyle}

\newtheorem{theorem}{Theorem}
\newtheorem{lemma}[theorem]{Lemma}
\newtheorem{prop}[theorem]{Proposition}
\newtheorem{cor}[theorem]{Corollary}
\numberwithin{equation}{section}
\numberwithin{theorem}{section}

\title{Wintenberger's Functor for Abelian Extensions}
\author{Kevin Keating \\
Department of Mathematics \\
University of Florida \\
Gainesville, FL 32611 \\
USA \\[.2cm]
{\tt keating@math.ufl.edu}}
\begin{document}

\maketitle

\begin{abstract}
Let $k$ be a finite field.  Wintenberger used the field of
norms to give an equivalence between a category whose objects
are totally ramified abelian $p$-adic Lie extensions $E/F$,
where $F$ is a local field with residue field $k$, and a
category whose objects are pairs $(K,A)$, where $K\cong k((T))$
and $A$ is an abelian $p$-adic Lie subgroup of $\Aut_k(K)$.
In this paper we extend this equivalence to allow
$\Gal(E/F)$ and $A$ to be arbitrary abelian pro-$p$ groups.
\end{abstract}

\section{Introduction}

     Let $q=p^f$ and let $k$ be a finite field with $q$
elements.  We define a category $\A$ whose objects are
totally ramified abelian extensions $E/F$, where $F$ is a
local field with residue field $k$, and $[E:F]$ is infinite if
$F$ has characteristic 0.  An $\A$-morphism from $E/F$ to
$E'/F'$ is defined to be a continuous embedding
$\rho:E\ra E'$ such that
\begin{enumerate}
\item $\rho$ induces the identity on $k$.
\item $E'$ is a finite separable extension of $\rho(E)$.
\item $F'$ is a finite separable extension of $\rho(F)$.
\end{enumerate}
Let $\rho^*:\Gal(E'/F')\ra\Gal(E/F)$ be the map induced by
$\rho$.  It follows from conditions 2 and 3 that $\rho^*$ has
finite kernel and finite cokernel.

     We also define a category $\B$ whose
objects are pairs $(K,A)$, where $K$ is a local field of
characteristic $p$ with residue field $k$ and $A$
is a closed abelian subgroup of $\Aut_k(K)$.  A
$\B$-morphism from $(K,A)$ to $(K',A')$ is defined to be
a continuous embedding $\sigma:K\ra K'$ such that
\begin{enumerate}
\item $\sigma$ induces the identity on $k$.
\item $K'$ is a finite separable extension of $\sigma(K)$.
\item $A'$ stabilizes $\sigma(K)$, and the image of $A'$ in
$\Aut_k(K)$ is an open subgroup of $A$.
\end{enumerate}
Let $\sigma^*:A'\ra A$ be the map induced by
$\sigma$.  It follows from conditions 2 and 3 that $\sigma^*$ has
finite kernel and finite cokernel.

     The field of norms construction \cite{cn} gives a
functor $\f:\A\ra\B$ defined by
\begin{equation} \label{functor}
\f(E/F)=(X_F(E),\Gal(E/F)).
\end{equation}
We wish to prove the following:

\begin{theorem} \label{main}
$\f$ is an equivalence of categories.
\end{theorem}

     Wintenberger (\cite{WZp,Wab}; see also \cite{lie}) has
shown that $\f$ induces
an equivalence between the full subcategory $\A_L$ of $\A$
consisting of extensions $E/F$ such that $\Gal(E/F)$ is an
abelian $p$-adic Lie group, and the full subcategory $\B_L$
of $\B$ consisting of pairs $(K,A)$ such that $A$ is an
abelian $p$-adic Lie group.  (Contrary to \cite{lie,Wab} we
consider finite groups to be $p$-adic Lie groups.  The equivalence
of categories proved in \cite{lie,WZp,Wab} extends trivially
to include the case of finite
groups.)  The proof of Theorem~\ref{main} is based on reducing
to the equivalence between $\A_L$ and $\B_L$.

     The following result, proved by Laubie \cite{ram},
is a consequence of Theorem~\ref{main}:

\begin{cor} \label{laubie}
Let $(K,A)\in\B$.  Then there is $E/F\in\A$ such that $A$ is
isomorphic to $G=\Gal(E/F)$ as a filtered group.  That is,
there exists an isomorphism $i:A\ra G$ such that
$i(A[x])=G[x]$ for all $x\ge0$, where $A[x]$, $G[x]$
denote the ramification subgroups of $A$, $G$ with respect to
the lower numbering.
\end{cor}

     The finite field $k\cong{\mathbb F}_q$ is fixed
throughout the paper, as is the field $K=k((T))$ of formal
Laurent series over $k$.  We work with complete discretely
valued fields $F$ whose residue field is identified with $k$,
and with totally ramified abelian extensions of such fields.
The ring of integers of $F$ is denoted by $\O_F$ and the maximal
ideal of $\O_F$ is denoted by $\M_F$.  We let $v_F$
denote the valuation on the separable closure
$F^{sep}$ of $F$ which is normalized so that
$v_F(F^{\times})=\Z$, and we let $v_p$ denote the $p$-adic
valuation on $\Z$.  We say that the profinite group $G$
is finitely generated if there is a finite set $S\subset G$ such
that $\langle S\rangle$ is dense in $G$.

\section{Ramification theory and the field of norms}

     In this section we recall some facts from ramification
theory, and summarize the construction of the field of norms
for extensions in $\A$.

     Let $E/F\in\A$.  Then $G=\Gal(E/F)$ has a decreasing
filtration by the upper ramification subgroups $G(x)$, defined
for nonnegative real $x$.  (See for instance \cite[IV]{cl}.)
We say that $u$ is an upper ramification break of $G$ if
$G(u+\epsilon)\subsetneqq G(u)$ for every $\epsilon>0$.
Since $G$ is abelian, by the Hasse-Arf Theorem
\cite[V\,\S7, Th.\,1]{cl} every upper ramification break
of $G$ is an integer.  In addition, since $E/F$ is a totally
ramified abelian extension, it follows from class field theory
that $E/F$ is arithmetically profinite (APF) in the sense of
\cite[\S1]{cn}.  This
means that for every $x\ge0$ the upper ramification subgroup
$G(x)$ has finite index in $G=G(0)$.  This allows us to
define the Hasse-Herbrand functions
\begin{equation}
\psi_{E/F}(x)=\int_0^x|G(0):G(t)|\,dt
\end{equation}
and $\phi_{E/F}(x)=\psi_{E/F}^{-1}(x)$.  It follows that the
ramification subgroups of $G$ with the lower numbering can be
defined by $G[x]=G(\phi_{E/F}(x))$ for $x\ge0$.
We say that $l$ is a lower
ramification break for $G$ if $G[u+\epsilon]\subsetneqq G[u]$
for every $\epsilon>0$.  It is clear from the definitions that
$l$ is a lower ramification break if and only if $\phi_{E/F}(l)$
is an upper ramification break.

     When $(K,A)\in\B$ the abelian subgroup $A$ of
$\Aut_k(K)$ also has a ramification filtration.
The lower ramification subgroups of $A$ are defined by
\begin{equation}
A[x]=\{\sigma\in A: v_K(\sigma(T)-T)\ge x+1\}
\end{equation}
for $x\ge0$.  Since $A[x]$ has finite index in $A=A[0]$ for
every $x\ge0$, the function
\begin{equation}
\phi_A(x)=\int_0^x\frac{dt}{|A[0]:A[t]|}
\end{equation}
is strictly increasing.  Hence we can define the ramification
subgroups of $A$ with the upper numbering by
$A(x)=A[\psi_A(x)]$, where $\psi_A(x)=\phi_A^{-1}(x)$.
The upper and lower ramification breaks
of $A$ are defined in the same way as the upper and lower
ramification breaks of $\Gal(E/F)$.  The lower ramification
breaks of $A$ are certainly integers, and Laubie's result
(Corollary~\ref{laubie}) together with the Hasse-Arf theorem
imply that the upper ramification breaks of $A$ are also integers.

     For $E/F\in\A$ let $i(E/F)$ denote the smallest (upper
or lower) ramification break of the extension $E/F$.
The following basic result from ramification theory is
presumably well-known (cf.~\cite[3.2.5.5]{cn}).

\begin{lemma} \label{breaks}
Let $M/F\in\A$, let $L\in\E_{M/F}$, and let $L'/L$ be a
finite totally ramified abelian extension which is linearly
disjoint from $M/L$.  Assume that $M'=ML'$ has residue field $k$,
so that $M'/F'\in\A$.  Then
$i(M'/F')\le\psi_{F'/F}(i(M/F))$, with equality if the largest
upper ramification break $u$ of $F'/F$ is less than $i(M/F)$.
\end{lemma}

\proof Set $G=\Gal(M'/F)$, $H=\Gal(M'/M)$, and
$N=\Gal(M'/F')$.  Then $G=HN\cong H\times N$.  Let
$y=\phi_{F'/F}(i(M'/F'))$.  Then
\begin{equation}
N=N(i(M'/F'))=N(\psi_{F'/F}(y))=G(y)\cap N.
\end{equation}
It follows that $G(y)\supset N$, and hence that
$G/H=G(y)H/H=(G/H)(y)$.  Therefore $y\le i(M/F)$, which implies
$i(M'/F')\le\psi_{F'/F}(i(M/F))$.

     If $u<i(M/F)$ then the group
\begin{equation}
(G/N)(i(M/F))=G(i(M/F))N/N
\end{equation}
is trivial.  It follows that $G(i(M/F))\subset N$, and hence
that
\begin{equation}
N(\psi_{F'/F}(i(M/F)))=G(i(M/F))\cap N=G(i(M/F)).
\end{equation}
The restriction map from $\Gal(M'/F')=N$ to $\Gal(M/F)\cong
G/H$ carries $G(i(M/F))$ onto
\begin{equation}
G(i(M/F))H/H=(G/H)(i(M/F))=G/H.
\end{equation}
Thus $N(\psi_{F'/F}(i(M/F))=N$, so we have
$i(M'/F')\ge\psi_{F'/F}(i(M/F))$.  Combining this with the
inequality proved above we get
$i(M'/F')=\psi_{F'/F}(i(M/F))$.
\qed \medskip

     Let $E/F\in\A$.  Since $E/F$ is an APF extension, the
field of norms of $E/F$ is defined:  Let $\E_{E/F}$ denote
the set of finite subextensions of $E/F$, and for
$L',L\in\E_{E/F}$ such that $L'\supset L$ let
$\n_{L'/L}:L'\ra L$ denote the norm map.  The field of norms
$X_F(E)$ of $E/F$ is defined to be the inverse limit of
$L\in\E_{E/F}$ with respect to the norms.  We denote an
element of $X_F(E)$ by $\alpha_{E/F}=(\alpha_L)_{L\in\E_{E/F}}$.
Multiplication in $X_F(E)$ is defined componentwise, and
addition is defined by the rule
$\alpha_{E/F}+\beta_{E/F}=\gamma_{E/F}$, where
\begin{equation}
\gamma_L=\lim_{L'\in\E_{E/L}}\n_{L'/L}(\alpha_{L'}+\beta_{L'})
\end{equation}
for $L\in\E_{E/F}$.  We embed $k$ into $X_F(E)$ as follows:
Let $F_0/F$ be the maximum tamely ramified subextension of
$E/F$, and for $\zeta\in k$ let $\tilde{\zeta}_{F_0}$ be the
Teichm\"uller lift of $\zeta$ in $\O_{F_0}$.  Note that for any
$L\in\E_{E/F_0}$ the degree of the extension $L/F_0$ is a power
of $p$.  Therefore there is a unique $\tilde{\zeta}_L\in L$ such
that $\tilde{\zeta}_L^{[L:F_0]}=\tilde{\zeta}_{F_0}$.
Define $f_{E/F}(\zeta)$ to be
the unique element of $X_F(E)$ whose $L$ component is
$\tilde{\zeta}_L$ for every $L\in\E_{E/F_0}$.  Then the map
$f_{E/F}:k\ra X_F(E)$ is a field embedding.  By choosing
a uniformizer for $X_F(E)$ we get a $k$-isomorphism
$X_F(E)\cong k((T))$.

     The ring of integers $\O_{X_F(E)}$ consists of those
$\alpha_{E/F}\in X_F(E)$ such that $\alpha_L\in\O_L$ for all
$L\in\E_{E/F}$ (or equivalently, for any $L\in\E_{E/F}$).  A
uniformizer $\pi_{E/F}=(\pi_L)_{L\in\E_{E/F}}$ for $X_F(E)$
gives a uniformizer $\pi_L$ for each finite subextension $L/F$
of $E/F$, and also a
uniformizer $\pi_{M/F}=(\pi_L)_{L\in\E_{M/F}}$ of $X_F(M)$
for each infinite subextension $M/F$ of $E/F$.
The action of $\Gal(E/F)$ on the
fields $L\in\E_{E/F}$ induces a $k$-action of $\Gal(E/F)$ on
$X_F(E)$.  By identifying $\Gal(E/F)$ with the subgroup of
$\Aut_k(X_F(E))$ which it induces, we get the functor
$\f(E/F)=(X_F(E),\Gal(E/F))$ which was mentioned in
(\ref{functor}).

     Let $E'$ be a finite extension of $E$ such that $E'/F\in\A$.
Then there is $M\in\E_{E/F}$ and a finite extension $M'$ of
$M$ such that $E'=EM'$ and $E$, $M'$ are linearly disjoint
over $M$.  We define an embedding $j:X_F(E)\ra X_F(E')$ as
follows.  For $\alpha_{E/F}\in X_F(E)$ set
$j(\alpha_{E/F})=\beta_{E'/F}$, where $\beta_{E'/F}$ is the
unique element of $X_F(E')$ such that $\beta_{LM'}=\alpha_L$
for all $L\in\E_{E/M}$ \cite[3.1.1]{cn}.  The
embedding $j$ makes $X_F(E')$ into a finite separable
extension of $X_F(E)$ of degree $[E':E]$; in this setting we
denote $X_F(E')$ by $X_{E/F}(E')$.
If $E''\supset E'\supset E$ are finite extensions such that
$E''/F\in\A$ then $X_{E/F}(E')/X_F(E)$ is a subextension
of $X_{E/F}(E'')/X_F(E)$.  Let $D/E$ be an infinite
extension such that $D/F\in\A$.  Then $X_{E/F}(D)$ is defined
to be the union of $X_{E/F}(E')$ as $E'$ ranges over the finite
subextensions of $D/E$.

     For $E/F\in\A$ set
$r(E/F)=\lceil\frac{p-1}{p}\cdot i(E/F)\rceil$.  The proof of
Theorem~\ref{main} depends on the following two results, the
first of which was proved by Wintenberger:

\begin{prop} \label{quotient}
Let $E/F\in\A$, let $L\in\E_{E/F_0}$, and define
$\xi_L:\O_{X_F(E)}\ra\O_L/\M_L^{r(E/L)}$ by
$\xi_L(\alpha_{E/F})=\alpha_L\pmod{\M_L^{r(E/L)}}$.  Then
\medskip \\
(a) $\xi_L$ is a surjective ring homomorphism. \smallskip \\
(b) $\xi_L$ induces the map $\zeta\mapsto\zeta^{p^{-a}}$ on $k$,
where $a=v_p([L:F])$.
\end{prop}

\proof This is proved in
Propositions 2.2.1 and 2.3.1 of \cite{cn}. \qed

\begin{prop} \label{commute}
Let $E/F\in\A$, let $L\in\E_{E/F}$, and let $L'/L$ be a
finite totally ramified abelian extension which is linearly
disjoint from $E/L$.  Assume that $E'=EL'$ has residue field $k$,
so that $E'/L'\in\A$.
Then the following diagram commutes, where the bottom horizontal
map is induced by the inclusion $\O_L\hookrightarrow \O_{L'}$:
\begin{equation} \label{diagram}
\begin{array}{ccc}
\O_{X_F(E)}&\overset{j}{\longrightarrow}&\O_{X_{E/F}(E')} \\[.3cm]
\!\!\!\!\!\!\xi_L\downarrow&&\downarrow\xi_{L'}\!\!\!\!\!\! \\[.3cm]
\O_L/\M_L^{r(E/L)}&\longrightarrow&
\O_{L'}/\M_{L'}^{r(E'/L')}
\end{array}
\end{equation}
Furthermore, for $\zeta\in k$ we have
$j\circ f_{E/F}(\zeta)=f_{E'/F}(\zeta^{p^b})$, where
$b=v_p([L':L])$.
\end{prop}

\proof Using Lemma~\ref{breaks} we get
\begin{equation}
i(E'/L')\le\psi_{L'/L}(i(E/L))\le[L':L]i(E/L).
\end{equation}
Thus $r(E'/L')\le[L':L]r(E/L)$, so the bottom horizontal map
in the diagram is well-defined.
Let $\alpha_{E/F}=(\alpha_M)_{M\in\E_{E/F}}$ be an element
of $\O_{X_F(E)}$.  Then $j(\alpha_{E/F})$ is the unique
element of $\O_{X_{E/F}(E')}$ whose $ML'$-component is
$\alpha_M$ for every $M\in\E_{E/L}$.  In particular, the
$L'$-component of $j(\alpha_{E/F})$ is $\alpha_L$.  Hence
$\xi_L(\alpha_{E/F})$ and $\xi_{L'}(j(\alpha_{E/F}))$ are both
congruent to $\alpha_L$ modulo $\M_{L'}^{r(E'/L')}$, which proves
the commutativity of (\ref{diagram}).  The second claim follows
from (\ref{diagram}) and Proposition~\ref{quotient}(b). \qed

\section{Proof of Theorem \ref{main}} \label{proof}

     In this section we prove Theorem~\ref{main}.  To do
this, we must show that the functor $\f:\A\ra\B$ is
essentially surjective and fully faithful.

     We begin by showing that $\f$ is essentially surjective.
Let $A$ be a closed abelian subgroup of $\Aut_k(K)$, where
$K=k((T))$.  Then $A$ is a $p$-adic Lie group if and only if
$A$ is finitely generated.  Since $\f$ induces an equivalence
between the categories $\A_L$ and $\B_L$, it suffices to prove
that $(K,A)$ lies in the essential image of
$\f$ in the case where $A$ is {\em not} finitely generated.

     Let $F\cong k((T))$, let $E/F$ be a finite totally
ramified abelian extension, and let $\pi$ be a uniformizer of
$E$.  Then for each $\sigma\in\Gal(E/F)$ there is a unique
$f_{\sigma}\in k[[T]]$ such that $\sigma(\pi)=f_{\sigma}(\pi)$.
Let $a=v_p([E:F])$ and define
\begin{equation}
G(E/F,\pi)=\{\gamma\in\Aut_k(K):\gamma(T)=f_{\sigma}^{p^a}(T)
\mbox{ for some }\sigma\in\Gal(E/F)\},
\end{equation}
where $f_{\sigma}^{p^a}(T)$ is the power series obtained
from $f_{\sigma}(T)$ by replacing the coefficients by their
$p^a$ powers.
Then $G(E/F,\pi)$ is a subgroup of $\Aut_k(K)$ which is
isomorphic to $\Gal(E/F)$.

     Let $l_0<l_1<l_2<\ldots$ denote the positive lower
ramification breaks of $A$.  For $n\ge0$ set
$r_n=\lceil\frac{p-1}{p}\cdot l_n\rceil$ and let
$\Gammabar_n$ denote the quotient of $\Gamma=\Aut_k(K)$ by
the lower ramification subgroup
\begin{equation}
\Gamma[r_n-1]=
\{\sigma\in\Gamma:\sigma(T)\equiv T\pmod{\M_K^{r_n}}\}.
\end{equation}
Let $\s_n$ denote the set of pairs $(E,\pi)$ such that
\begin{enumerate}
\item $E/F$ is a totally ramified abelian subextension of
$F^{sep}/F$ such that $\Gal(E/F)[l_n]$ is trivial.
(Such an extension is necessarily finite.)
\item $\pi$ is a uniformizer of $E$ such that the image of
$G(E/F,\pi)$ in $\Gammabar_n$ is equal to the image of $A$ in
$\Gammabar_n$.
\end{enumerate}
Since there are only finitely many extensions $E/F$
satisfying condition 1, and condition 2 depends only on the
class of $\pi$ modulo $\M_E^{r_n}$, the set $\s_n$ is
compact.  Using the following lemma we get a map from $\s_n$
to $\s_{n-1}$.

\begin{lemma} \label{nu}
Let $n\ge1$, let $(E,\pi)\in\s_n$, and let $\tilde{E}$ denote
the fixed field of $\Gal(E/F)[l_{n-1}]$.  Then
$(\tilde{E},\n_{E/\tilde{E}}(\pi))\in\s_{n-1}$.
\end{lemma}

\proof It follows from the definitions that $\tilde{E}/F$
is a totally ramified abelian extension and that
$\Gal(\tilde{E}/F)[l_{n-1}]$ is trivial.
Set $\tilde{\pi}=\n_{E/\tilde{E}}(\pi)$, choose
$\sigma\in\Gal(E/F)$, and let $\tilde{\sigma}$ denote the
restriction of $\sigma$ to $\tilde{E}$.  By
\cite[Prop.\,2.2.1]{cn} the norm $\n_{E/\tilde{E}}$ induces
a ring homomorphism from $\O_E$ to
$\O_{\tilde{E}}/\M_{\tilde{E}}^{r_{n-1}}$.  Therefore
\begin{equation} \label{sigtilde}
\tilde{\sigma}(\tilde{\pi})=\n_{E/\tilde{E}}(\sigma(\pi))=
\n_{E/\tilde{E}}(f_{\sigma}(\pi))\equiv
f_{\sigma}^{p^b}(\n_{E/\tilde{E}}(\pi))\pmod{\M_{\tilde{E}}^{r_{n-1}}},
\end{equation}
where $b=v_p([E:\tilde{E}])$.  Let $\tilde{a}=v_p([\tilde{E}:F])$
and let $f_{\tilde{\sigma}}\in k[[T]]$ be such that
$\tilde{\sigma}(\tilde{\pi})=f_{\tilde{\sigma}}(\tilde{\pi})$.
Then by (\ref{sigtilde}) we have
\begin{alignat}{2}
f_{\tilde{\sigma}}(T)&\equiv
f_{\sigma}^{p^{b}}(T)&&\pmod{T^{r_{n-1}}} \\
f_{\tilde{\sigma}}^{p^{\tilde{a}}}(T)&\equiv
f_{\sigma}^{p^{a}}(T)&&\pmod{T^{r_{n-1}}}.
\end{alignat}
It follows that
$G(\tilde{E}/F,\tilde{\pi})$ and $G(E/F,\pi)$ have the same
image in $\Gammabar_{n-1}$, and hence that
$G(\tilde{E}/F,\tilde{\pi})$ and $A$ have the same image in
$\Gammabar_{n-1}$.
\qed \medskip

     Since each $A/A[l_n]$ is finite there is a sequence
$A_0\le A_1\le A_2\le\dots$ of finitely generated closed
subgroups of $A$ such that $A_nA[l_n]=A$ for all $n\ge0$.
Recall that $\f$ induces an equivalence of categories between
$\A_L$ and $\B_L$.  Since $(K,A_n)\in\B_L$, for $n\ge0$ there
exists $L_n/F_n\in\A_L$ such that $\f(L_n/F_n)$ is
$\B$-isomorphic to $(K,A_n)$.  Since $A_n$ is a normal
subgroup of $A$,
the action of $A$ on $K$ gives a $\B$-action of $A$ on the
pair $(K,A_n)$.  Since $\f(L_n/F_n)\cong(K,A_n)$ and $\f$ induces
an equivalence between $\A_L$ and $\B_L$, this action is induced
by a faithful $\A$-action of $A$ on $L_n/F_n$.  Since
$\Gal(L_n/F_n)\cong A_n$ is finitely generated, and $A$ is not
finitely generated, this implies that $\Aut_k(F_n)$ is not
finitely generated.  Hence $F_n$ has characteristic $p$.
Therefore we may assume $F_n=F$ and $L_n\subset F^{sep}$
with $F\cong k((T))$ fixed.

     For $n\ge0$ let
\begin{equation}
i_n:(K,A_n)\lra(X_{F}(L_n),\Gal(L_n/F))
\end{equation}
be an $\A$-isomorphism, and set $\pi_{L_n/F}=i_n(T)$.  Let
$E_n\subset L_n$ be the fixed field of $\Gal(L_n/F)[l_n]$.
Then either $E_n=L_n$, or $i(L_n/E_n)\ge l_n$, in which case
$r(L_n/E_n)\ge r_n$.  Therefore by
Proposition~\ref{quotient}(a) we have
$(E_n,\pi_{E_n})\in\s_n$, so $\s_n\not=\varnothing$.
For $n\ge1$ let $\nu_n:\s_n\ra\s_{n-1}$ be the map defined
by Lemma~\ref{nu}.  Since each $\s_n$ is compact, by
Tychonoff's theorem there exists a sequence of pairs
$(E_n,\pi_{E_n})\in\s_n$ such that
$\nu_n(E_n,\pi_{E_n})=(E_{n-1},\pi_{E_{n-1}})$ for $n\ge1$.
It follows in particular that
$F\subset E_0\subset E_1\subset E_2\subset\dots$.
Let $E_{\infty}=\cup_{n\ge 0}E_n$.  Then $E_{\infty}$ is
a totally ramified abelian extension of $F$, and the
uniformizers $\pi_{E_n}$ for $E_n$ induce a
uniformizer $\pi_{E_{\infty}/F}$ for $X_F(E_{\infty})$.
Let $\tau$ denote the unique $k$-isomorphism from $K=k((T))$
to $X_F(E_{\infty})$ such that $\tau(T)=\pi_{E_{\infty}/F}$.
It follows from our construction that $\tau$ induces a
$\B$-isomorphism from $(K,A)$ to
\begin{equation}
\f(E_{\infty}/F)=(X_F(E_{\infty}),\Gal(E_{\infty}/F)).
\end{equation}
Thus $(K,A)$ lies in the essential image of $\f$, so $\f$ is
essentially surjective.

     We now show that $\f$ is faithful.  Let $E/F$ and
$E'/F'$ be elements of $\A$, and set $G=\Gal(E/F)$ and
$G'=\Gal(E'/F')$.  We need to show that the map
\begin{equation} \label{Hom}
\Psi:\Hom_{\A}(E/F,E'/F')\lra
\Hom_{\B}((X_F(E),G),(X_{F'}(E'),G'))
\end{equation}
induced by the field of norms functor is one-to-one.
Suppose $\rho_1,\rho_2\in\Hom_{\A}(E/F,E'/F')$ satisfy
$\Psi(\rho_1)=\Psi(\rho_2)$.  Let
$\pi_{E/F}=(\pi_L)_{L\in\E_{E/F}}$ be a uniformizer for $X_F(E)$.
Then $\Psi(\rho_1)(\pi_{E/F})=\Psi(\rho_2)(\pi_{E/F})$, and hence
$(\rho_1(\pi_L))_{L\in\E_{E/F}}=(\rho_2(\pi_L))_{L\in\E_{E/F}}$.
It follows that $\rho_1(\pi_L)=\rho_2(\pi_L)$ for every
$L\in\E_{E/F}$.  Since $\rho_1$ and $\rho_2$ are $k$-algebra
homomorphisms, this implies that $\rho_1=\rho_2$.

     It remains to show that $\f$ is full, i.\,e., that
$\Psi$ is onto.  It follows from the arguments given in
the proof of \cite[Th.\,2.1]{Wab} that the
codomain of $\Psi$ is empty if $\char(F)\not=\char(F')$,
and that $\Psi$ is onto if $G$ and $G'$
are finitely generated.  In particular, $\Psi$
is onto if $\char(F)=0$ or $\char(F')=0$.
If one of $G$, $G'$ is finitely generated
and the other is not then the codomain of $\Psi$ is
empty.  Hence it suffices to prove that $\Psi$
is onto in the case where $\char(F)=\char(F')=p$ and neither of
$G$, $G'$ is finitely generated.

     We first show that every isomorphism lies in the image of
$\Psi$.  Let
\begin{equation}
\tau:(X_F(E),G)\lra(X_{F'}(E'),G')
\end{equation}
be a $\B$-isomorphism.  For $n\ge1$ let $F_n$ denote the
fixed field of $G[l_n]=G(u_n)$.  If
$\dst\lim_{n\ra\infty}l_n/[F_n:F]=\infty$ then an argument
similar to that used in \cite[\S2]{WZp} shows
that $\tau$ is induced by an $\A$-isomorphism from $E/F$ to
$E'/F'$.  This limit condition holds for instance if $\char(F)=p$
and $\Gal(E/F)$ is finitely generated, but it can fail
if $\Gal(E/F)$ is not finitely generated.
Therefore we use a different method to prove that
$\tau$ lies in the image of $\Psi$, based on a
characterization of $F_n/F$ in terms of $(X_F(E),G)$.

     Let $d$ denote the $F_n$-valuation of the different of
$F_n/F$, and let $c$ be an integer such that
$c>\phi_{F_n/F}(\frac{p}{p-1}(l_{n-1}+d))$.  Since $G/G(c)$ is
finite there exists a finitely generated closed
subgroup $H$ of $G$ such that $HG(c)=G$.  Let
$M\subset E$ be the fixed field of $H$ and set $M_n=F_nM$.  Then
$F_n/F$ and $M_n/M$ are finite abelian extensions.  On the other
hand, since $G$ is not finitely generated, $\Gal(M/F)\cong G/H$
is not finitely generated, and hence $M/F$ is an infinite abelian
extension.

\begin{prop} \label{kisom}
Let $\pi_{E/F}$ be a uniformizer for $X_F(E)$, and recall that
$\pi_{E/F}$ induces uniformizers $\pi_F$, $\pi_{F_n}$,
$\pi_{M/F}$, and $\pi_{M_n/F}$ for $F$, $F_n$, $X_F(M)$, and
$X_{M/F}(M_n)$.  There exists a $k$-isomorphism
$\zeta:X_{M/F}(M_n)/X_F(M)\ra F_n/F$ such that
\begin{enumerate}
\item $\zeta(\pi_{M/F})=\pi_F$;
\item $\zeta(\pi_{M_n/F})\equiv\pi_{F_n}
\pmod{\M_{F_n}^{l_{n-1}+1}}$;
\item $\gamma\cdot\zeta(\pi_{M_n/F})=\zeta(\gamma\cdot\pi_{M_n/F})$
for every $\gamma\in H$.
\end{enumerate}
\end{prop}

     The proof of this proposition depends on the following
lemma (cf.~\cite[p.\,88]{trunc}).

\begin{lemma} \label{near}
Let $F$ be a local field, let $g(T)\in\O_F[T]$ be a
separable monic
Eisenstein polynomial, and let $\alpha\in F^{sep}$ be a root of
$g(T)$.  Set $E=F(\alpha)$ and let $d=v_E(g'(\alpha))$ be the
$E$-valuation of the different of the extension $E/F$.  Choose
$\pi\in F^{sep}$ such that $v_E(g(\pi))>d$.  Then there is a root
$\beta$ for $g(X)$ such that $v_E(\pi-\beta)\ge v_E(g(\pi))-d$.
\end{lemma}

\proof Let $\alpha_1,\alpha_2,\dots,\alpha_n$ be the
roots of $g(T)$, and choose $1\leq j\leq n$ to maximize
$w=v_E(\pi-\alpha_j)$.  For $1\leq i\le n$ we have
\begin{equation}
v_E(\pi-\alpha_i) \geq \min\{w,v_E(\alpha_j-\alpha_i)\},
\end{equation}
with equality if $w>v_E(\alpha_j-\alpha_i)$.  Since $w\geq
v_E(\pi-\alpha_i)$, this implies that for $i\not=j$ we have
$v_E(\pi-\alpha_i)\leq v_E(\alpha_j-\alpha_i)$.  Since
\begin{equation}
g(\pi)=(\pi-\alpha_1)(\pi-\alpha_2)\dots(\pi-\alpha_n),
\end{equation}
we get
\begin{equation} \label{m}
v_E(g(\pi))\leq
w+\sum_{\substack{1\le i\le n\\[.1cm]i\not=j}}
\,v_E(\alpha_j-\alpha_i)=w+d.
\end{equation}
Setting $\beta=\alpha_j$ gives $v_E(\pi-\beta)=w\ge v_E(g(\pi))-d$.
\qed \medskip

\noindent
{\em Proof of Proposition \ref{kisom}:}
Since $HG(c)=G$ we have
\begin{equation}
i(M/F)\ge c>\phi_{F_n/F}(l_{n-1})=u_{n-1}.
\end{equation}
Therefore by Lemma~\ref{breaks} we get
\begin{equation}
i(M_n/F_n)=\psi_{F_n/F}(i(M/F))\ge\psi_{F_n/F}(c),
\end{equation}
and hence $r(M_n/F_n)\ge s$, where $s=\lceil\frac{p-1}{p}
\cdot\psi_{F_n/F}(c)\rceil$.  Let $g(T)$ be the
minimum polynomial for $\pi_{M_n/F}$ over $X_F(M)$,
and let $g_F(T)\in\O_F[T]$ be the
polynomial obtained by applying the canonical map
$\xi_F:X_F(M)\ra F$ to the coefficients of $g(T)$.  Since
$g(\pi_{M_n/F})=0$, it follows from Propositions~\ref{quotient}(a)
and \ref{commute} that $v_{F_n}(g_F(\pi_{F_n}))\ge
r(M_n/F_n)\ge s$.  On the other
hand, let $\mu:X_F(M)\ra F$ be the unique $k$-algebra
isomorphism such that $\mu(\pi_{M/F})=\pi_F$.  Then by
Proposition~\ref{quotient}(a) we have
$\mu(\alpha_{M/F})\equiv\alpha_F\pmod{\M_F^t}$ for all
$\alpha_{M/F}\in \O_{X_F(M)}$, where
$t=\lceil\frac{p-1}{p}\cdot c\rceil$.
Let $g^{\mu}(T)\in\O_F[T]$ be the polynomial obtained by applying
$\mu$ to the coefficients of $g(T)$.  Then
$g^{\mu}(T)\equiv g_F(T)\pmod{\M_F^t}$.  Since $[F_n:F]\cdot c\ge
\psi_{F_n/F}(c)$ we have $[F_n:F]\cdot t\ge s$; therefore
since $v_{F_n}(g_F(\pi_{F_n}))\ge s$ we have
$v_{F_n}(g^{\mu}(\pi_{F_n}))\ge s>l_{n-1}+d$.
It follows from Lemma~\ref{near}
that there is a root $\beta$ of $g^{\mu}(T)$ such that
$v_{F_n}(\pi_{F_n}-\beta)>l_{n-1}$.  Therefore by Krasner's
Lemma we have $F(\beta)\supset F(\pi_{F_n})$.  Since
\begin{equation}
[F(\beta):F]=\deg(g)=[F(\pi_{F_n}):F]
\end{equation}
we deduce that
$F(\beta)=F(\pi_{F_n})=F_n$.  Since $\pi_{M_n/F}$ is a root
of $g(T)$, and $\beta$ is a root of $g^{\mu}(T)$, the isomorphism
$\mu$ from $X_F(M)$ to $F$ extends uniquely to an isomorphism
$\zeta$ from $X_{M/F}(M_n)/X_F(M)$ to $F_n/F$ such that
$\zeta(\pi_{M_n/F})=\beta\equiv\pi_{F_n}
\pmod{\M_{F_n}^{l_{n-1}+1}}$.

     We now show that $\zeta$ is $H$-equivariant.
Let $\gamma\in H$ and define $\psi_{\gamma}\in k[[T]]$ by
\begin{equation}
\psi_{\gamma}(\pi_{M_n/F})=\gamma\cdot\pi_{M_n/F}=
(\gamma\cdot\pi_L)_{L\in\E_{M_n/F}},
\end{equation}
where we identify $k$ with a subfield of
$X_F(M)$ using the map $f_{M/F}$.  Using
Propositions~\ref{quotient} and \ref{commute} we get
\begin{equation}
\gamma\cdot\pi_{F_n}\equiv \psi_{\gamma}(\pi_{F_n})
\pmod{\M_{F_n}^{r(M_n/F_n)}}.
\end{equation}
Since $\zeta(\pi_{M_n/F})\equiv\pi_{F_n}\pmod{\M_{F_n}^{l_{n-1}+1}}$
and $r(M_n/F_n)\ge s\ge l_{n-1}+1$ this implies
\begin{alignat}{2}
\zeta(\gamma\cdot\pi_{M_n/F})&=\zeta(\psi_{\gamma}(\pi_{M_n/F})) \\
&=\psi_{\gamma}(\zeta(\pi_{M_n/F})) \\
&\equiv \psi_{\gamma}(\pi_{F_n})&&\pmod{\M_{F_n}^{l_{n-1}+1}} \\
&\equiv \gamma\cdot\pi_{F_n}&&\pmod{\M_{F_n}^{l_{n-1}+1}} \\
&\equiv\gamma\cdot\zeta(\pi_{M_n/F})&&
\pmod{\M_{F_n}^{l_{n-1}+1}}.
\end{alignat}
Since $\zeta(\gamma\cdot\pi_{M_n/F})$ and
$\gamma\cdot\zeta(\pi_{M_n/F})$ are both roots of $g^{\mu}(T)$
we deduce that $\gamma\cdot\zeta(\pi_{M_n/F})=
\zeta(\gamma\cdot\pi_{M_n/F})$.  Since $\zeta$ is $k$-linear
and $\gamma$ acts trivially on $k$, it follows that
$\gamma\cdot\zeta(\alpha)=\zeta(\gamma\cdot\alpha)$ for all
$\alpha\in X_{M/F}(M_n)$.~\qed \medskip

     Since $\tau$ is an $\A$-isomorphism, $\tau^*:G'\ra G$ is
a group isomorphism.  For $\gamma\in G$ set
$\gamma'=(\tau^*)^{-1}(\gamma)$, and for $N\le G$ set
$N'=(\tau^*)^{-1}(N)$.  Then $\tau$ induces an isomorphism
from $(X_F(E),N)$ to $(X_{F'}(E'),N')$.  In particular,
$\tau$ gives an isomorphism from $(X_F(E),H)$ to
$(X_{F'}(E'),H')$.  Using the isomorphism
$X_{X_F(M)}(X_{M/F}(E))\cong X_F(E)$ from \cite[3.4.1]{cn} we
get an isomorphism
\begin{equation} \label{tauH}
\tau_H:(X_{X_F(M)}(X_{M/F}(E)),H)\lra
(X_{X_{F'}(M')}(X_{M'/F'}(E')),H'),
\end{equation}
where $M'\subset E'$ is the fixed field of $H'$.  Since
$H$ is an abelian $p$-adic Lie group, it follows from
\cite{lie,WZp,Wab} that $\tau_H$ is induced by an
$\A$-isomorphism
\begin{equation}
\rho:X_{M/F}(E)/X_F(M)\lra X_{M'/F'}(E')/X_{F'}(M').
\end{equation}
By restricting $\rho$ we get an isomorphism
\begin{equation}
\tilde{\rho}:X_{M/F}(M_n)/X_F(M)\lra X_{M'/F'}(M_n')/X_{F'}(M'),
\end{equation}
where $M_n'=(M')_n=F_n'M'$ is the fixed field of
$H'[l_n]=H[l_n]'$.  Furthermore, for
$\gamma\in H$ and $\alpha\in X_{M/F}(M_n)$ we have
$\tilde{\rho}(\gamma(\alpha))=\gamma'(\tilde{\rho}(\alpha))$.

     Let $\pi_{E/F}$ be a uniformizer for $X_F(E)$, set
$\pi_{E'/F'}=\tau(\pi_{E/F})$, and let
\begin{align}
\zeta&:X_{M/F}(M_n)/X_F(M)\lra F_n/F \\
\zeta'&:X_{M'/F'}(M_n')/X_{F'}(M')\lra F_n'/F'
\end{align}
be the isomorphisms given by Proposition~\ref{kisom}.  Then
$\omega_n=\zeta'\circ\tilde{\rho}\circ\zeta^{-1}$ gives a
$k$-linear isomorphism from $F_n/F$ to $F_n'/F'$.  It follows
from Proposition~\ref{kisom} that
\begin{equation} \label{scong}
\omega_n(\pi_{F_n})\equiv\pi_{F_n'}\pmod{\M_{F_n'}^{l_{n-1}+1}},
\end{equation}
and that for all $\gamma\in H$ we have
\begin{equation} \label{sequi}
\omega_n(\gamma(\pi_{F_n}))=\gamma'(\omega_n(\pi_{F_n})).
\end{equation}
Since the restriction map from $H=\Gal(E/M)$ to $\Gal(F_n/F)$ is
onto, (\ref{sequi}) is actually valid for all $\gamma\in G$.

     Let $\T_n$ denote the set of $k$-isomorphisms
$\omega_n:F_n/F\ra F_n'/F'$ which satisfy (\ref{scong})
and (\ref{sequi}) for all $\gamma\in G$.  Since $l_{n-1}$
is the only ramification break of $F_n'/F_{n-1}'$ we have
$\psi_{F_n'/F_{n-1}'}(l_{n-1})=l_{n-1}$.  Therefore by
(\ref{scong}) and \cite[V\,\S6, Prop.\,8]{cl}, for any
$\omega_n\in\T_n$ we have
\begin{equation} \label{norm}
\n_{L_n'/L_{n-1}'}(\omega_n(\pi_{F_n}))\equiv
\n_{L_n'/L_{n-1}'}(\pi_{F_n'})\pmod{\M_{F_{n-1}'}^{l_{n-1}+1}}.
\end{equation}
Since $\n_{L_n/L_{n-1}}(\pi_{F_n})=\pi_{F_{n-1}}$ and
$\n_{L_n'/L_{n-1}'}(\pi_{F_n'})=\pi_{F_{n-1}'}$, it follows
from (\ref{norm}) and (\ref{sequi}) that
\begin{equation}
\omega_n(\pi_{F_{n-1}})\equiv\pi_{F_{n-1}'}
\pmod{\M_{F_{n-1}'}^{l_{n-1}+1}}.
\end{equation}
Therefore restriction induces a map from $\T_n$ to $\T_{n-1}$.

     Define a metric on $\T_n$ by setting
$d(\omega_n,\tilde{\omega}_n)=2^{-a}$, where
$a=v_{F_n'}(\omega_n(\pi_{F_n})-\tilde{\omega}_n(\pi_{F_n}))$.
Then $\T_n$ is easily seen to be compact, and we showed above
that $\T_n$ is nonempty.  Therefore by Tychonoff's theorem
there is a sequence $(\omega_n)_{n\ge1}$ with
$\omega_n\in\T_n$ and $\omega_n|_{F_{n-1}}=\omega_{n-1}$.
Since $E=\cup_{n\ge1}F_n$ and $E'=\cup_{n\ge1}F_n'$ the
isomorphisms $\omega_n:F_n/F\ra F_n'/F'$ combine to give
an $\A$-isomorphism $\Omega:E/F\ra E'/F'$.  Let
$\theta=\Psi(\Omega)$
be the $\B$-isomorphism induced by $\Omega$ and let
$m_n=\min\{l_{n-1}+1,r(E/F_n)\}$.  It follows from
(\ref{scong}) and Proposition~\ref{quotient}(a) that
\begin{equation}
\theta(\pi_{E/F})\equiv\pi_{E'/F'}\pmod{\M_{X_{F'}(E')}^{m_n}}
\end{equation}
for every $n\ge1$.  Since $\lim_{n\ra\infty}m_n=\infty$
we get $\theta(\pi_{E/F})=\pi_{E'/F'}=\tau(\pi_{E/F})$.
Hence $\tau=\theta=\Psi(\Omega)$.

     Now let $\sigma$ be an arbitrary element of
$\Hom_{\B}((X_F(E),G),(X_{F'}(E'),G'))$.  Since $X_{F'}(E')$
is a finite separable extension of $\sigma(X_F(E))$, by
\cite[3.2.2]{cn} there is a finite separable extension
$\tilde{E}'$ of $E$ such that $\sigma$ extends to an
isomorphism
\begin{equation}
\tau:X_{E/F}(\tilde{E}')\lra X_{F'}(E').
\end{equation}
It follows that each $\gamma'\in G'$ induces an automorphism
$\tilde{\gamma}$ of $X_{E/F}(\tilde{E}')$ whose restriction to
$X_F(E)$ is $\sigma^*(\gamma')\in G$.  Since $X_{E/F}(F^{sep})$
is a separable closure of $X_F(E)$ \cite[Cor.\,3.2.3]{cn},
$\tilde{\gamma}$ can be extended to an automorphism $\gammabar$
of $X_{E/F}(F^{sep})$.  Since $\gammabar$ stabilizes $X_F(E)$,
and $\gammabar|_{X_F(E)}=\sigma^*(\gamma')$ is induced by an
element of $G=\Gal(E/F)$, it follows from
\cite[Rem.\,3.2.4]{cn} that $\gammabar$ is induced by an
element of $\Gal(F^{sep}/F)$, which we also denote by
$\gammabar$.  Since $\gammabar$ stabilizes
$X_{E/F}(\tilde{E}')$, it stabilizes $\tilde{E}'$ as well.
Thus $\gammabar|_{\tilde{E}'}$ is an element of
$\Aut_k(\tilde{E}')$ which is uniquely determined by
$\gamma'$.  Since $\gammabar|_{\tilde{E}'}$ induces the
automorphism $\tilde{\gamma}$ of
$X_{E/F}(\tilde{E}')$, we denote $\gammabar|_{\tilde{E}'}$
by $\tilde{\gamma}$ as well.

     Let $\tilde{F}'$ denote the subfield of $\tilde{E}'$
which is fixed by the group
$\tilde{G}'=\{\tilde{\gamma}:\gamma'\in G'\}$.
Then $\tilde{F}'\supset F$, so $\tilde{E}'/\tilde{F}'$ is a
Galois extension.  Since the image of $\tilde{G}'\cong G'$
in $G$ is open,
and $\tilde{E}'/E$ is a finite extension, it follows that
$\tilde{F}'$ is a finite separable extension of $F$, and
$\tilde{G}'=\Gal(\tilde{E}'/\tilde{F}')$.  In particular,
$\tilde{F}'\cong k((T))$ is a local field with residue
field $k$.  Hence $(X_{\tilde{F}'}(\tilde{E}'),\tilde{G}')$
is an object in $\B$, and $\tau$ gives a $\B$-isomorphism
from $(X_{\tilde{F}'}(\tilde{E}'),\tilde{G}')$ to
$(X_{F'}(E'),G')$.  By the arguments given above
we see that $\tau$ is induced by an $\A$-isomorphism
$\Omega:\tilde{E}'/\tilde{F}'\ra E'/F'$.
Furthermore, the embedding $E\hookrightarrow\tilde{E}'$
induces an $\A$-morphism $i:E/F\ra\tilde{E}'/\tilde{F}'$.  Let
\begin{equation}
\alpha:(X_F(E),G)\lra(X_{\tilde{F}'}(\tilde{E}'),\tilde{G}')
\end{equation}
be the $\B$-morphism induced by $i$.  Then
$\sigma=\tau\circ\alpha=\Psi(\Omega\circ i)$.


\begin{thebibliography}{9}

\bibitem{trunc} K. Keating, Extensions of local fields and
truncated power series, J. Number Theory {\bf 116} (2006),
69--101. 

\bibitem{lie} F. Laubie, Extensions de Lie et groupes
d'automorphismes de corps locaux, Compositio Math. {\bf 67}
(1988), 165--189.

\bibitem{ram} F. Laubie, Ramification des groupes ab\'eliens
d'automorphismes de $\F_q((X))$, Canad.\ Math.\ Bull.\ 
{\bf 50} (2007), 594--597.

\bibitem{cl} J.-P. Serre, {\em Corps Locaux}, Hermann, Paris
(1962).

\bibitem{WZp} J.-P.~Wintenberger, Automorphismes des corps
locaux de charact\'eristique $p$, J. Th\'eor.\ Nombres
Bordeaux {\bf 16} (2004), 429--456.

\bibitem{Wab} J.-P.~Wintenberger, Extensions ab\'eliennes et
groupes d'automorphismes de corps locaux,
C. R. Acad.\ Sci.\ Paris S\'er. A-B {\bf290} (1980),
A201--A203.

\bibitem{cn} J.-P.~Wintenberger,  Le corps des normes de
certaines extensions infinies de corps locaux; applications,
Ann.\ Sci.\ \'Ecole Norm.\ Sup.\ (4) {\bf 16} (1983), 59--89.

\end{thebibliography}
\end{document}